\documentclass[12pt,a4paper,reqno]{amsart}
\usepackage{amsmath}
\usepackage{amsthm}
\usepackage{amsfonts}
\usepackage{amssymb}
\usepackage{setspace}
\numberwithin{equation}{section}

\theoremstyle{plain}

\newtheorem{lemma}[subsection]{Lemma}

\theoremstyle{definition}

\renewcommand{\leq}{\leqslant}
\renewcommand{\geq}{\geqslant}

\newsavebox{\proofbox}
\savebox{\proofbox}{\begin{picture}(7,7)%
 \put(0,0){\framebox(7,7){}}\end{picture}}

\newcommand{\md}[1]{\ensuremath{(\mbox{mod}\, #1)}}

\def\E{\mathbb{E}}
\def\Z{\mathbb{Z}}
\def\R{\mathbb{R}}
\def\T{\mathbb{T}}

\DeclareMathOperator{\vol}{vol}

\def\ud{\,\mathrm{d}}

\def\proof{\noindent\textit{Proof. }}
\def\endproof{\hfill{\usebox{\proofbox}}}

\begin{document}

\onehalfspace
\title{A note on Elkin's improvement of Behrend's construction}

\author{Ben Green}
\address{Centre for Mathematical Sciences\\
Wilberforce Road\\
     Cambridge CB3 0WA\\
     England
}
\email{b.j.green@dpmms.cam.ac.uk}

\author{Julia Wolf}
\address{Mathematical Sciences Research Institute\\
17 Gauss Way\\
Berkeley, CA 94720, U.S.A.
}
\email{julia.wolf@cantab.net}

\thanks{The first author holds a Leverhulme Prize and is grateful to the Leverhulme Trust for their support. This paper was written while the authors were attending the special semester in ergodic theory and additive combinatorics at MSRI}

\begin{abstract}
We provide a short proof of a recent result of Elkin in which large subsets of $\{1,\dots,N\}$ free of 3-term progressions are constructed.
\end{abstract}
\maketitle

\begin{center}
\emph{To Mel Nathanson}
\end{center}

\section{Introduction}

Write $r_3(N)$ for the cardinality of the largest subset of $\{1,\dots,N\}$ not containing three distinct elements in arithmetic progression. A famous construction of Behrend \cite{behrend} shows, when analysed carefully, that 
\[ r_3(N) \gg \frac{1}{\log^{1/4} N} \cdot \frac{N}{2^{2\sqrt{2}\sqrt{\log_2 N}}}.\]
In a recent preprint \cite{elkin} Elkin was able to improve this 62-year old bound to
\[ r_3(N) \gg \log^{1/4} N \cdot \frac{N}{2^{2\sqrt{2}\sqrt{\log_2 N}} }.\]
Our aim in this note is to provide a short proof of Elkin's result. It should be noted that the only advantage of our approach is brevity: it is based on ideas morally close to those of Elkin, and moreover his argument is more constructive than ours.

Throughout the paper $0 < c < 1 < C$ denote absolute constants which may vary from line to line. We write $\T^d = \R^d/\Z^d$ for the $d$-dimensional torus.

\section{The proof}

Let $d$ be an integer to be determined later, and let $\delta \in (0,1/10)$ be a small parameter (we will have $d  \sim C\sqrt{\log N}$ and $\delta \sim \exp{(-C\sqrt{\log N}}$)). Given $\theta,\alpha \in \T^d$, write $\Psi_{\theta,\alpha} : \{1,\dots,N\} \rightarrow \T^d$ for the map $n \mapsto \theta n + \alpha\md{1}$.

\begin{lemma}\label{lemma1}
Suppose that $n$ is an integer. Then $\Psi_{\theta,\alpha}(n)$ is uniformly distributed on $\T^d$ as $\theta,\alpha$ vary uniformly and independently over $\T^d$. Moreover, if $n$ and $n'$ are distinct positive integers, then the pair $(\Psi_{\theta,\alpha}(n),\Psi_{\theta,\alpha}(n'))$ is uniformly distributed on $\T^d \times \T^d$ as $\theta,\alpha$ vary uniformly and independently over $\T^d$. 
\end{lemma}
\proof Only the second statement requires an argument to be given. Perhaps the easiest proof is via Fourier analysis, noting that 
\[ \int e^{2\pi i(k \cdot (\theta n + \alpha) + k'\cdot (\theta n' + \alpha))} \ud \theta \ud \alpha = 0\]
unless $k + k' = kn + k'n' = 0$. Provided that $k$ and $k'$ are not both zero, this cannot happen for distinct positive integers $n,n'$. Since the exponentials $e^{2\pi i(kx + k'x')}$ are dense in $L^2(\T^d \times \T^d)$, the result follows.\endproof\vspace{11pt}

Let us identify $\T^d$ with $[0,1)^d$ in the obvious way. For each $r \leq \frac{1}{2}\sqrt{d}$, write  $S(r)$ for the region
\[ \{ x\in [0,1/2]^d : r - \delta \leq \Vert x \Vert_2\leq r\}.\]
\begin{lemma}\label{vol-lemma} There is some choice of $r$ for which $\vol(S(r)) \geq c\delta 2^{-d}$.
\end{lemma}
\proof First note that if $(x_1,\dots,x_d)$ is chosen at random from $[0,1/2]^d$ then, with probability at least $c$, we have $ | \Vert x \Vert_2 - \sqrt{d/12}| \leq C$. This is a consequence of standard tail estimates for sums of independent identically distributed random variables, of which $\Vert x \Vert_2^2 = \sum_{i=1}^d x_i^2$ is an example. The statement of the lemma then immediately follows from the pigeonhole principle.\endproof\vspace{11pt}

Write $S := S(r)$ for the choice of $r$ whose existence is guaranteed by the preceding lemma; thus $\vol(S)  \geq c \delta 2^{-d}$. Write $\tilde S$ for the same set $S$ but considered now as a subset of $[0,1/2]^d \subseteq \R^d$.  Since there is no ``wraparound'', the 3-term progressions in $S$ and $\tilde{S}$ coincide and henceforth we abuse notation, regarding $S$ as a subset of $\R^d$ and dropping the tildes. (To use the additive combinatorics jargon, $S$ and $\tilde{S}$ are \emph{Freiman isomorphic}.) Suppose that $(x,y) $ is a pair for which $x-y,x$ and $x+y$ lie in $S$. By the parallelogram law
\[ 2\Vert x \Vert_2^2 + 2\Vert y \Vert_2^2 = \Vert x + y \Vert_2^2 + \Vert x - y \Vert_2^2\] and straightforward algebra we have
\[ \Vert y \Vert_2 \leq \sqrt{r^2 - (r-\delta)^2} \leq \sqrt{2\delta r}.\]
It follows from the formula for the volume of a sphere in $\R^d$ that the volume of the set $B \subseteq \T^d \times \T^d$ in which each such pair $(x,y)$ must lie is at most $\vol(S) C^d (\delta/\sqrt{d})^{d/2}$.

The next lemma is an easy observation based on Lemma \ref{lemma1}.
\begin{lemma}\label{lemma3} Suppose that $N$ is even.
Define $A_{\theta,\alpha} := \{n \in [N] : \Psi_{\theta,\alpha}(n) \in S\}$. Then
\[ \E_{\theta,\alpha} |A_{\theta,\alpha}| =  N\vol(S) \]
whilst the expected number of nontrivial 3-term arithmetic progressions in $A_{\theta,\alpha}$ is 
\[ \E_{\theta,\alpha} T(A_{\theta,\alpha}) = \frac{1}{4}N(N-5)\vol(B).\]
\end{lemma}
\proof The first statement is an immediate consequence of the first part of Lemma \ref{lemma1}.
Now each nontrivial 3-term progression is of the form $(n- d,n,n + d)$ with $d \neq 0$. Since $N$ is even there are $N(N-5)/4$ choices for $n$ and $d$, and each of the consequent progressions lies inside $A_{\theta,\alpha}$ with probability $\vol(B)$ by the second part of Lemma \ref{lemma1}.\endproof\vspace{11pt}

To finish the argument, we just have to choose parameters so that 
\begin{equation}\label{to-satisfy} \frac{1}{3}\vol(S) \geq \frac{1}{4}(N-5)\vol(B).\end{equation}Then we shall have 
\[ \E \big(\frac{2}{3} |A_{\theta,\alpha}| - T(A_{\theta,\alpha}) \big)\geq \frac{1}{3}N\vol(S).\]
In particular there is a specific choice of $A := A_{\theta,\alpha}$ for which both $T(A) \leq 2|A|/3$ and $|A| \geq \frac{1}{2}N\vol(S)$. Deleting up to two thirds of the elements of $A$, we are left with a set of size at least $\frac{1}{6}N\vol(S)$ that is free of 3-term arithmetic progressions.

To do this it suffices to have $ C^d(\delta/\sqrt{d})^{d/2} \leq c/N$, which can certainly be achieved by taking $\delta := c\sqrt{d} N^{-2/d}$. For this choice of parameters we have, by the earlier lower bound on $\vol(S)$, that
\[ |A| \geq \frac{1}{6}N\vol(S) \geq c\sqrt{d}2^{-d} N^{1-2/d}.\] Choosing $d := \lceil\sqrt{2\log_2 N}\rceil$ we recover Elkin's bound.\endproof

\section{A question of Graham}

The authors did not set out to try and find a simpler proof of Elkin's result. Rather, our concern was with a question of Ron Graham (personal communication to the first-named author, see also \cite{graham,landman}). Defining $W(2;3,k)$ to be the smallest $N$ such that any red-blue colouring of $[N]$ contains either a 3-term red progression or a $k$-term blue progression, Graham asked whether $W(2;3;k) < k^A$ for some absolute constant $A$ or, even more ambitiously, whether $W(2;3,k) \leq Ck^2$. Our initial feeling was that the answer was surely no, and that a counterexample might be found by modifying the Behrend example in such a way that its complement does not contain long progressions. Reinterpreting the Behrend construction in the way that we have done here, it seems reasonably clear that it is not possible to provide a negative answer to Graham's question in this way. 

\section{acknowledgement}

The authors are grateful to Tom Sanders for helpful conversations.

\providecommand{\bysame}{\leavevmode\hbox to3em{\hrulefill}\thinspace}
\providecommand{\MR}{\relax\ifhmode\unskip\space\fi MR }
\providecommand{\MRhref}[2]{%
  \href{http://www.ams.org/mathscinet-getitem?mr=#1}{#2}
}
\providecommand{\href}[2]{#2}

     \end{document}